%% file: v2.tex
\documentclass[12pt]{amsart}
\usepackage[margin=1in]{geometry}
\usepackage[utf8]{inputenc}
\usepackage{amsmath}
\usepackage{amssymb}
\usepackage{esint}

\usepackage[pdftitle={Plateau Problem},
  pdfauthor={Samuel Bronstein},
  pdffitwindow=true,
  breaklinks=true,
  colorlinks=true,
  urlcolor=blue,
  linkcolor=red,
  citecolor=red,
  anchorcolor=red]{hyperref}
\usepackage{titleps}
\usepackage{xcolor}
\definecolor{links}{RGB}{70,0,255}
\hypersetup{citecolor=blue,linkcolor=links}

\usepackage{subfiles}

\numberwithin{equation}{section}
\numberwithin{table}{section}
\numberwithin{figure}{section}

\newcommand{\C}{\mathbb C}

\newcommand{\R}{\mathbb R}
\renewcommand{\H}{\mathbb H}
\renewcommand{\S}{\mathbb S}
\newcommand{\X}{\mathbb X}

\newcommand{\Tr}{\text{Tr}}

\newcommand{\calN}{\mathcal N}

\newcommand{\eps}{\varepsilon}

\newcommand{\II}{\text{I\!I}}

\newcommand{\PSL}{\mathrm{PSL}}

\theoremstyle{plain}
\newtheorem{theorem}{Theorem}[section]
\newtheorem*{theorem*}{Theorem}
\newtheorem{lemma}[theorem]{Lemma}
\newtheorem{proposition}[theorem]{Proposition}

\newtheorem{corollary}[theorem]{Corollary}

\theoremstyle{definition}
\newtheorem{definition}[theorem]{Definition}

\newtheorem{problem}[theorem]{Problem}

\theoremstyle{remark}
\newtheorem{remark}{Remark}[section]

\title{On almost-fuchsian submanifolds of Hadamard spaces and the asymptotic Plateau problem}
\author{samuel bronstein}
\date{2023}

\begin{document}
\begin{abstract}
We consider minimal submanifolds of negatively curved spaces with small curvature.
We show that in a Hadamard space with negatively pinched curvature $-C\leq K\leq -1$,
complete minimal submanifolds with second fundamental form less than $1$ everywhere 
bound a class of spheres at infinity for which the asymptotic
Plateau problem is uniquely solvable.
\end{abstract}
\maketitle
\tableofcontents
\section{Introduction}
\subfile{v2intro}
\section{The convex hull barrier}
\subfile{CVHull}
\section{Almost-fuchsian submanifolds}
\subfile{v2AFsub}
\section{The asymptotic Plateau problem}
\subfile{v2uniq}

\newpage
\bibliographystyle{alpha}
\bibliography{references}
\end{document}

%% file: v2intro.tex
In 1983, Uhlenbeck~\cite{Uhl83} introduced the notion of \emph{almost-fuchsian representation}.
A representation of a surface group in~$\PSL(2,\C)$ is almost-fuchsian if it is discrete,
faithful and there is an equivariant minimal disc in $\H^3$ whose principal values lie
in a compact set of $(-1,1)$.
Since then, almost-fuchsian representations have been a very fruitful field of study,
see for instance \cite{Eps86, KS07, HW13, Sep16, San17, Jia21, ES22, BS23}
We consider here a generalization of almost-fuchsian submanifolds in negatively curved spaces:
Let $\X$ be a complete simply connected riemannian manifold whose sectional curvature is less
or equal than $-1$.
\begin{definition}
An immersion $f:Y\rightarrow\X$ is \emph{almost-fuchsian} if it is minimal,
proper, and it satisfies
\begin{equation}
\sup_{y\in Y}|\II_f|<1
\end{equation}
\end{definition}
Remark that $Y$ can always be endowed with the induced metric which
will be complete by the properness assumption.
By extension, a submanifold $Y\subset\X$ will be said to be almost-fuchsian if the inclusion
map is almost-fuchsian.

In this paper, we prove the following embedding theorem:
\begin{theorem}
Let $f:Y\rightarrow\X$ be an almost-fuchsian immersion.
Then $Y$ is contractible. Denote by $NY$ the normal bundle to $Y$.
The exponential map $\exp:NY\rightarrow\X$ is a diffeomorphism.
Also, $f$ extends at infinity to an embedding from a sphere $S\rightarrow\partial_\infty\X$.
\end{theorem}

We also consider the following asymptotic Plateau problem: given a
sphere~$S\subset\partial_\infty\X$, can we count the complete minimal submanifolds of $\X$
asymptotically bounding $S$ ? Jiang~\cite{Jia21} and Huang--Lowe--Seppi~\cite{HLS23}
considered this problem for spheres in $\partial_\infty\H^n$ bounding respectively
an almost-fuchsian disc or an almost-fuchsian hypersurface. They proved that in that case,
the unique minimal submanifold of the given dimension bounding $S$ is the almost-fuchsian one,
which is unique.

Here we prove that this statement holds in a broader generality, provided that $\X$
has negatively pinched curvature between $-C$ and $-1$.

\begin{theorem}
Let $\X$ be an $n$-dimensional negatively pinched Hadamard space with sectional curvature
less or equal than~$-1$. Let $Y\subset\X$ be a $k$-dimensional almost-fuchsian submanifold
bounding a $(k-1)$-dimensional sphere $S\subset\partial_\infty\X$.
Then $Y$ is the unique $k$-dimensional complete minimal submanifold of $\X$
asymptotically bounding $S$.
\end{theorem}

We don't know whether this statement holds when the sectional curvature of $\X$ has no lower
bound, because we lack of convex sets to work with.
Note that Huang--Lowe--Seppi \cite{HLS23} also cover the case of weakly-almost-fuchsian
discs in $\H^3$, that is when the supremum of the principal values is allowed to be one.

This paper is divided in three parts. The first part is about a very standard fact, that
minimal submanifolds of a space $\X$ with negatively pinched sectional curvature
remain in the convex hull of their asymptotic boundary. Surprisingly, we couldn't prove it 
without a lower bound on the curvature of $\X$.
The second part is devoted to the proof of the embedding theorem, with some
explicit estimates on the geometry of an almost-fuchsian submanifold.
Finally, the third part is devoted to the proof of the uniqueness to the asymptotic Plateau
problem for spheres bounding an almost-fuchsian submanifold.

The author thanks Andrea Seppi, Graham Smith and Nicolas Tholozan for their help, insights and
discussions on these topics.

%% file: CVHull.tex
This section is devoted to the proof of the following theorem:
\begin{theorem}\label{CvxHullIsBarrier}
Let $\X$ be a Hadamard space with negatively pinched curvature.
Let $Y$ be a minimal submanifold of $\X$, bounding a subset $F\subset\partial_\infty\X$.
Denote $C$ the convex hull of $F$.
Then $Y\subset C$.
\end{theorem}
While this theorem is well known when $\X=\H^n$, we didn't find a general proof of this fact in
negative curvature.
A very elegant proof in $\H^n$ uses totally geodesic hypersurfaces as barriers, as the convex
hull is the intersection of all half spaces containing $C$ \cite{Jia21}.
While this proof can be word for word done in a symmetric space of a semisimple Lie group of
rank~$1$, here we need to be more careful, as totally geodesic hypersurfaces might not exist.

In all this section $\X$ will denote a negatively pinched Hadamard space,
its boundary $\partial_\infty\X$ will be its ideal boundary which equals its Gromov boundary
in that case \cite{BGS85,Gro78,EO73}.
The main idea is to use the following property:
\begin{proposition}\label{CvxRegion}
Let $F$ be a closed set of $\partial_\infty\X$, and denote by $C$ its convex hull.
Let $p\in X-C$.
Then there is a closed domain $D\subset\X\cup\partial_\infty\X$ disjoint from
$C\cup F$ such that $p\in D$, and $D$ is foliated by smooth strictly convex hypersurfaces.
\end{proposition}
\begin{proof}
We use a theorem of Bowditch \cite{Bow94}. As $\X$ has negatively pinched curvature,
there is a metric on $\X\cup\partial_\infty\X$ compatible with the topologies
such that the map $F\mapsto CH(F)$ is continuous on the closed sets of $\X\cup\partial_\infty\X$
equipped with the Hausdorff distance.
Therefore, there is a neighborhood $V$ of $F$ such that $\widetilde C$ the convex hull
of $V$ does not contain $p$.
The complement of $\widetilde C$ is then foliated by the equidistant to $\widetilde C$,
which are strictly convex.

Now a priori the convex $\widetilde C$ has no smooth boundary, and so our foliation is
only $C^{1,1}$. However, Parkkonen--Paulin \cite{PP12} (Proposition 6) proved that
a $\eps$-neighborhood of $\widetilde C$ contains a smooth convex set, so up
to replacing $\widetilde C$ by this smooth convex set we assume the foliation of the complement
to be by smooth convex hypersurfaces. Denote $D$ the closure of the complement of $\widetilde C$.
\end{proof}
\begin{proof}[Proof of theorem \ref{CvxHullIsBarrier}]
Let $Y$ be a minimal submanifold and let $C$ denote the convex hull of $\partial_\infty Y$.
Let $p\in\X-C$ and $D$ denote the domain provided by proposition \ref{CvxRegion}.

Suppose that $Y$ intersects $D$.
Because $D$ does not bound $\partial_\infty Y$, the intersection $Y\cap D$ is compact.
Then there must some point $q$ at which $Y$ is tangent to one of the strictly convex hypersurfaces
foliating $D$. This contradicts the minimality of $Y$. Hence $Y\cap D$
is empty.
Repeating this argument for any $p\in\X-C$ proves that $Y\subset C$.
\end{proof}
\begin{remark}
The lower bound assumption on the curvature of $\X$ is crucial here,
otherwise there might not be enough convex sets to pursue our proof.
See for instance Ancona \cite{Anc94}, there are spaces with curvature less than $-1$ such that
the convex hull of any nontrivial open set at infinity is the full space $\X$.
\end{remark}

%% file: v2AFsub.tex
This part is about the notion of almost-fuchsian submanifold. First introduced for
surface group representations in~$\PSL(2,\C)$ by Uhlenbeck \cite{Uhl83}, almost-fuchsian discs
and hypersurfaces have a rich literature \cite{Eps86,KS93,Sep16,San17,Jia21,ES22,BS23,HLS23}.

Throughout this section~$\X$ will be a Hadamard space with sectional curvature less or equal
than~$-1$. Note that we do not require its curvature to admit a lower bound.
We consider an immersion~$f:Y\rightarrow\X$, where $Y$ is connected.
We assume $f$ to be complete, i.e. the induced metric on $Y$ is complete.

We will prove the following theorem.
\begin{theorem}\label{AFemb}
Assume $f:Y\rightarrow\X$ is minimal and satisfies:
\begin{equation}
\sup_{y\in Y} |\II_f|<1
\end{equation}
Then $Y$ is contractible. Denote by $NY$ the normal bundle to $Y$.
The exponential map $\exp:NY\rightarrow\X$ is a diffeomorphism.
Also, $f$ extends at infinity to an embedding from a sphere $S\rightarrow\partial_\infty\X$
\end{theorem}
The proof is divided in several steps.
First, we show that the induced metric by $\exp$ is nondegenerate on $NY$.
This will imply that $\exp$ is a covering map to $\X$, and so it is a diffeomorphism.
Finally, we will see that $f$ is a quasi-isometric embedding from $Y$ to $\X$,
so it extends at infinity to an embedding $S\rightarrow\partial_\infty\X$.

First note that because of the condition on $f$, the induced metric has negatively pinched
curvature on $Y$, and it makes sense to talk about its boundary $\partial_\infty Y$, which is
constructed in the works of Gromov \cite{Gro78}, Eberlein--O'NEill \cite{EO73} and
 Ballmann--Gromov--Schroeder \cite{BGS85}.
\begin{lemma}\label{NegPinch}
Let $f:Y\rightarrow\X$ be an almost-fuchsian immersion.
Let~$\eps>0$ such that
\begin{equation}
\sup_{y\in Y}|\II_f|\leq 1-\eps\,.
\end{equation}
Endow $Y$ with the induced metric by $f$.
Then its sectional curvature is less than $-\eps(2-\eps)$.
If the sectional curvature of $\X$ is lower bounded by $-C$,
then the sectional curvature of $f^\ast g_\X$ is bigger than $-C-2(1-\eps)^2$.
\end{lemma}
\begin{remark}
When $Y$ is a surface, the minimality condition ensures that we have
more precise bounds:
The sectional curvature of an almost-fuchsian surface is then less than $-1$
and bigger than~$-C-(1-\eps)^2$.
\end{remark}
\begin{proof}
This is a consequence of the Gauss condition:
For $u,v$ an orthonormal basis of a tangent plane to $Y$,
By the Gauss' equation,
the curvature of~$f^\ast g_\X$ is then
\begin{equation}
R_{f^\ast g_\X}(u,v,v,u)=R_{g_\X}(u,v,v,u)+g_\X(\II_f(u,u),\II_f(v,v))
-g_\X(\II_f(u,v),\II_f(u,v))
\end{equation}
As $R_{g_\X}(u,v,v,u)\leq -1$ and $|\II_f(u,u)|\leq 1-\eps$, we directly get
\begin{equation}
R_{f^\ast g_\X}(u,v,v,u)\leq -1+(1-\eps)^2=-\eps(2-\eps)\,.
\end{equation}
as desired.
Moreover, if $R_{g_\X}(u,v,v,u)\geq -C$,
we directly get that 
\begin{equation}
R_{f^\ast g_\X}(u,v,v,u)\geq -C-2(1-\eps)^2\,,
\end{equation}
as claimed.
\end{proof}
\begin{proposition}\label{CharEq}
Let $f:Y\rightarrow\X$ be an almost-fuchsian immersion.
Denote by $NY$ the normal bundle to $Y$, and note $G$ the induced metric by the exponential
map on $NY$. For $t>0$, denote by $N_tY$ the tangent bundle to the distance $t$ hypersurface to $Y$.
Let $(x, tv)\in NY $ with $v$ a unit normal vector.
Let $B_t$ be the  parallel transport map from $T_x\X\cap(v)^\perp$ to $T_{\exp{tv}}(N_tY)$.
Let $w\in T_v(UY)$ and let $f$ be the function:
\begin{equation}
f(t)=|B_tw|^2
\end{equation}
If $\X$ is the hyperbolic space, $f$ satisfies the following controls:
\begin{equation}
(\sqrt{f})''=\sqrt{f}
\end{equation}
If $\X$ has sectional curvature strictly less than $-1$, then
\begin{equation}
(\sqrt{f})''>\sqrt{f}
\end{equation}
Also, in both cases
\begin{equation}
(\ln f)''+\frac{(\ln f)'^2}{2}\geq 2
\end{equation}
with the same equality and inequality conditions.
\end{proposition}
\begin{proof}
With the introduced notations, fix $\gamma:I\rightarrow UY$ a smooth path such that
$\gamma(0)=v$ and $\gamma'(0)=w$. The family $c_s(t)$ is a smooth variation of
geodesics, hence $J(t)=\partial_sc_s(t)|_{s=0}$ is a Jacobi field along $c_0$.
Evaluate the characteristic equation of Jacobi fields against $J$:
\begin{equation}\label{Jaceq}
g_\X(\ddot J,J)+R_\X(J,\partial_t c,\partial_t c,J)=0
\end{equation}
In the second term, we recognize the sectional curvature of the tangent plane spanned
by $J$ and $\partial_t c$.
Choose a local orthonormal chart of $T_x\X$ such that $B_0=I_k\oplus O_{n-k-1}$ and 
$\dot B_0= B(v)\oplus I_{n-k-1}$, where $B(v)$ is the shape operator of $f$ at $v$.
Then equation \ref{Jaceq} rewrites as
\begin{equation}
^\perp w\ddot B_tB_t w+R_\X(J,\partial_t c,\partial_t c, J)=0
\end{equation}
and as the sectional curvature is less than $-1$, we replace the second term with the control
\begin{equation}
R_\X(J,\partial_t c,\partial_t c,J)\leq -^\perp w^\perp B_tB_t w
\end{equation}
to get the following
\begin{equation}
^\perp w ^\perp(\ddot B_t-B_t)B_t w\geq 0
\end{equation}
Now remark that
\begin{equation}
\left\{\begin{array}{ll}
f'(t)&=2^\perp w^\perp B_t \dot B_t w\\
f''(t)&=2^\perp w^\perp B_t \ddot B_t w +2^\perp w ^\perp \dot B_t\dot B_t w
\end{array}\right.
\end{equation}
So the characteristic equation rewrites as
\begin{equation}\label{FullEq}
f''(t)\geq 2f(t)+\frac{f'(t)^2}{2f(t)}
\end{equation}
with equality when $\X$ is the hyperbolic space, and strict inequality if $\X$ has sectional
curvature strictly less than $-1$.
As a direct consequence,
\begin{equation}
(\sqrt{f})''\geq\sqrt{f}
\end{equation}
with the same equality and inequality conditions.
Finally, introduce $g=\ln(f)$.
As long as it is well defined,
we deduce from \ref{FullEq} that
\begin{equation}
g''+\frac{(g')^2}{2}\geq 2
\end{equation}
\end{proof}
\begin{remark}
When $\X$ has varying curvature less or equal than $-1$,
up to rescaling the metric we can assume its sectional curvature to be strictly less than $-1$.
\end{remark}
\begin{proof}[Proof of theorem \ref{AFemb}]
With the notations introduced, Fix $w\in T_v(UY)$
and consider the function $f$ defined in proposition \ref{CharEq}.
Denote $w_1$ the projection of $w$ on $TY$, and assume it is nonzero.
The initial conditions, by construction of $f$,
\begin{equation}
\left\{
\begin{array}{ll}
f(0)&=\alpha^2=|w_1|^2>0\\
f'(0)&=2\alpha\beta=2^\perp w_1 B(v)w_1<2\alpha
\end{array}
\right.
\end{equation}
From proposition \ref{CharEq}, we deduce that
\begin{equation}
\left\{\begin{array}{lll}
f(t)&\geq (\alpha\cosh(t)+\beta\sinh(t))^2\\\label{fCont}
f'(t)&\geq 2(\alpha\cosh(t)+\beta\sinh(t))(\alpha\sinh(t)+\beta\cosh(t))\\
g'(t)&= \frac{f'(t)}{f(t)}\geq 2\frac{\alpha\sinh(t)+\beta\cosh(t)}{\alpha\cosh(t)+\beta\sinh(t)}
\end{array}\right.
\end{equation}
with equality when $\X$ is hyperbolic, and strict inequality when the curvature of
$\X$ is stricly less than $-1$.

A first consequence is that $f$ never vanishes. Hence $(B_t)$ has trivial kernel,
and the metric induced by $\exp$ on $NY$ is nondegenerate.
So $\exp$ is a local diffeomorphims from $NY$ to $\X$.

Moreover, because the immersion $Y\rightarrow\X$ is complete, it is proper.
Hence for any point of $x\in\X$ there is at least one point in $Y$ which minimizes the distance
to $x$. This implies that the map $\exp$ is surjective.
As a proper surjective local diffeomorphism,
$\exp$ satisfies the path lifting property, so it is a
covering map. But $Y$ is connected and $\X$ is contractible, so $\exp$ is a global diffeomorphism,
$Y$ is embedded and contractible too.

It remains to prove that $f$ is a quasi-isometric embedding.
In order to prove we introduce a metric on $NY$ for which the inverse map of the exponential
map will be Lipschitz, and such that $Y$ sits in $NY$ totally geodesically with the
induced metric by $f$:

Consider $(x,tv)\in NY$ and decompose its tangent space into
$T_{(x,tv)}NY=T_t(x,\R v)\oplus T_xY\oplus T_v\S^{n-k-1}$.

Introduce the metric $h$ on $NY$ which in that decomposition is:
\begin{equation}
h_{(x,tv)}=1\oplus \cosh^2(t)f^\ast g_{\X,x}\oplus \sinh^2(t)g_{\S^{n-k-1}}
\end{equation}
By construction, $Y$ sits inside $NY$ totally geodesically with the metric induced by $f$, that
we denote $g_Y$.

The equation on $f$ \ref{fCont} ensures that:
\begin{equation}
\exp^\ast g_\X\geq 1\oplus(\cosh(t)I_k+B(v)\sinh(t))^2g_Y\oplus\sinh^2(t)g_{\S^{n-k-1}}
\end{equation}
Denote $\delta=1-\sup|\II_f|>0$.
Then
\begin{equation}
\exp^\ast g_\X\geq \delta^2(1\oplus\cosh^2(t)g_Y\oplus\sinh^2(t)g_{\S^{n-k-1}})
=\delta^2 h
\end{equation}
Hence the inverse map of the exponential map is $\frac{1}{\delta}$-Lipschitz.

As a consequence, $f$ is a quasi-isometric embedding, as it satisfies:
\begin{equation}
d_\X(f(x_1),f(x_2))\leq d_Y(x_1,x_2)\leq\frac{1}{\delta}d_\X(f(x_1),f(x_2))
\end{equation}
Combined with the fact that $Y$ has negatively pinched sectional curvature proven in
lemma~\ref{NegPinch}, it ensures that $f$
extends at infinity to an embedding from a sphere $\S^{k-1}\rightarrow\partial_\infty\X$.
\end{proof}

As a corollary of the computations we made, we have some explicit estimates of the geometry
of $NY$, endowed with the induced metric by $\exp$.
\begin{corollary}\label{ExpCvx}
Let $f:Y\rightarrow\X$ be a complete almost-fuchsian immersion, with $\X$ a Hadamard
space.
Denote by $N_tY$ the distance $t$ hypersurface to $Y$.
At a point $\exp_x(tv)\in N_tY$, consider $\lambda_1\leq\ldots\leq\lambda_k$ the eigenvalues
of $B(v)$ the shape operator of $Y$ evaluated at $(x,v)$.

	Then the second fundamental form $\II_t$ of $N_tY$ has eigenvalues $\lambda_1^t\leq\ldots\leq\lambda_{n-1}^t$ which satisfy
\begin{flalign}
	\forall i\leq k,\quad \sum_{j=1}^i\lambda_j^t &\geq\sum_{j=1}^i
	\frac{\lambda_j+\tanh(t)}{1+\lambda_j\tanh(t)}\\
	\forall i> k,\quad \sum_{j=1}^i\lambda_j^t &\geq\sum_{j=1}^k
	\frac{\lambda_j+\tanh(t)}{1+\lambda_j\tanh(t)}+(i-k)\frac{1}{\tanh(t)}
\end{flalign}
In particular, $N_tY$ is $k$-convex in the sense of Sha \cite{Sha86}.
\end{corollary}
\begin{proof}
This is a consequence of the control on $\frac{f'}{f}$ gotten in \ref{fCont}.
The trace of $\II_t$ on a $k$-plane $P$
is the sum of $\frac{1}{2}\frac{f'(t)}{f(t)}$ for $w$ spanning a basis of the plane $P$.

The almost-fuchsian condition ensures that the sum of the $\lambda_i$ is 0,
and that they live in a compact set of $(-1,1)$, which ensures the $k$-convexity of
$\II_t$.
\end{proof}

%% file: v2uniq.tex
The Plateau problem was introduced in 1847, originally about a mathematical proof of the
existence of soap films bounding a wire frame. It was solved by Douglas \cite{Dou31} and
Rado~\cite{Rad30}.
Here we are interestes in a noncompact analogous problem: the asymptotic Plateau problem.
Fix $\X$ a Hadamard space of sectional curvature less than $-1$, the asymptotic Plateau problem
can be formulated as
\begin{problem}[Asymptotic Plateau Problem]
Let $S$ be a $(k-1)$-sphere embedded in $\partial_\infty\X$.
Is there a minimal $k$-dimensional submanifold $M$ of $\X$ asymptotically bounding $S$ ?
Is $M$ unique ?
\end{problem}
The asymptotic Plateau problem has been studied by Anderson~\cite{And83}. The problem
has since had a rich history, see the survey of Cozkunuzer~\cite{Cos09}.

Recently, Huang--Lowe--Seppi~\cite{HLS23} considered the asymptotic Plateau problem for a class
of almost-fuchsian submanifolds of $\H^n$. If a sphere in $\partial_\infty\H^n$ bounds an
almost-fuchsian hypersurface, then it is the unique minimal hypersurface bounding it.
Jiang~\cite{Jia21} also solved the asymptotic Plateau problem for Jordan curves bounding
an almost-fuchsian disc in $\partial_\infty\H^n$. In this section, we generalize these results
to any sphere bounding an almost-fuchsian submanifold in a space $\X$, under the assumption
that $\X$ has negatively pinched curvature.

\begin{theorem}[Asymptotic Plateau Problem for almost-fuchsian spheres]\label{PlaPbAF}
Let $\X$ be a negatively pinched Hadamard space with sectional curvature less or equal
than $-1$.
Let $Y\subset\X$ be an almost-fuchsian submanifold of dimension $k$, bounding at infinity
a $(k-1)$-dimensional sphere $S\subset\partial_\infty\X$.
Then $Y$ is the unique complete $k$-dimensional minimal submanifold of $\X$ bounding $S$.
\end{theorem}

In all the remainder of the article, $\X$ is assumed to have pinched sectional curvature
between $-C$ and $-1$. $Y\subset\X$ is a complete almost-fuchsian $k$-dimensional
submanifold, and $S$ is its asymptotic boundary in $\partial_\infty\X$.

The first step in the proof is that eventual other minimal submanifolds bounding $S$
remain at bounded distance from $Y$:
\begin{proposition}\label{BddDist}
Let $Z$ be a minimal submanifold of $\X$ bounding $S$.
Then $Z\subset \calN_r Y$, where $\calN_r Y$ is the $r$-uniform neighborhood of $Y$,
and $r=\tanh^{-1}(1-\sup|\II_f|)$.
\end{proposition}
\begin{proof}
Thanks to theorem \ref{CvxHullIsBarrier},
we know that $Z$ is included in the convex
hull of $S$.
But the explicit bounds shown in corollary \ref{ExpCvx} show that for
$r\geq\tanh^{-1}(1-\sup|\II_f|)$, the uniform neighborhoods $\calN_r(Y)$ are convex,
and contain $Y$ so bound $S$ too.
As a consequence, $\calN_r(Y)$ contains the convex hull of $Y$ which contains $Z$.
\end{proof}
\begin{remark}
The same corollary \ref{ExpCvx} shows that all uniform neighborhoods are $k$-convex,
which directly implies that the function $z\in Z\mapsto d(z,Y)$ cannot
have a local maximum. When the induced metric of $Z$ has lower bounded sectional curvature,
we could apply the Omori--Yau~\cite{Omo67,Yau75} maximum principle to conclude that $Z=Y$.
Here we don't prove that $Z$ has lower bounded sectional curvature, but
we prove that the maximum principle is still applicable. This property is called
\emph{stochastical completeness}, cf Pigola--Rigoli--Setti \cite{PRS05}.
\end{remark}
We want to apply a maximum principle to the distance function to $Y$ restricted on a minimal
submanifold $Z$, so we first prove it satisfies a strong subharmonicity condition.
\begin{proposition}\label{SupHarmFct}
Let $Z\subset\X$ be a minimal $k$-dimensional submanifold.
Consider the function $u:z\in Z\mapsto d(z,Y)^2$.
Then there is $C>0$ such that
\begin{equation}
\Delta u\geq C u
\end{equation}
\end{proposition}
\begin{proof}

	Consider $d$ the distance function to $Y$, defined on the whole space $X$.
	Because $Y$ is almost-fuchsian, for any point of $\X-Y$, the distance is attained
	at a unique point of $\X$, so the function $d$ is smooth on $X-Y$.

	Furthermore, at a point $x$ where $d(x)=t>0$,
	we have the formula
	\begin{equation}
		\nabla^2 d=0\oplus \II_{N_t(Y)}
	\end{equation}
	the $0$ is simply because $d$ is linear in the direction of the minimizing geodesic to $Y$.
	Using the formula
	\begin{equation}
		\nabla^2(d^2)=2\nabla d\cdot\nabla d+ 2d\nabla^2 d
	\end{equation}
	We deduce that
	\begin{equation}
		\nabla^2(u)\geq 2\oplus 2d\II_{N_t(Y)}
	\end{equation}
	Now, because $Z$ is a minimal submanifold of $X$, the laplacian of the restriction of $u$
	equals the trace of the restriction of its Hessian to the tangent space of $Z$
	\begin{equation}
		\Delta(u|_Z)=\Tr((\nabla^2 u)|_{TZ})
	\end{equation}
	Hence the bounds on the eigenvalues of $\II_{N_t(Y)}$ computed in corollary \ref{ExpCvx}
	and the control
	\begin{equation}
		\frac{\tanh(t)+\lambda_i}{1+\lambda_i\tanh(t)}\leq 1
	\end{equation}
	allow us to get the control:
	\begin{equation}
		\Delta(u|_Z)\geq 2 d\underset{\lambda_1+\ldots+\lambda_k=0,|\lambda_i|\leq 1-\eps}{\inf}
		\sum_{i=1,\ldots,k}\frac{\tanh(d)+\lambda_i}{1+\tanh(d)\lambda_i}=2d\Phi(d)
	\end{equation}
	Remark that $\Phi(d)$ is continuous, positive, vanishes only when $d=0$,
	and satisfies
	\begin{equation}
		\Phi'(0)\geq k(1-(1-\sup|\II_f|)^2)>0
	\end{equation}
	As $d$ is bounded by a constant $r=\tanh^{-1}(1-\sup|\II_f|)$,
	there is a constant $C>0$ such that on $[0,r]$, $\Phi(d)\geq C d$.
	We deduce
	that the restriction of $u$ to $Z$ satisfies
	\begin{equation}
		\Delta u\geq 2C u
	\end{equation}
	, as claimed.
\end{proof}
Now we want to prove that we can apply the maximum principle to $u$.
In order to do so, we will use the Khas'minskii test \cite{Kha60}, as stated in
Theorem~$3.1$ and Proposition~$3.2.$ of \cite{PRS05}.
\begin{theorem}\label{KhasMinskiiTest}
	Let $M,g$ be a Riemannian manifold, and assume that $M$ supports a $C^2$ function
	$\gamma$, which tends to infinity at infinity, and satisfies
	\begin{equation}
		\Delta\gamma\leq\lambda\gamma\text{ off a compact set}
	\end{equation}
	for some $\lambda>0$. Then $M,g$ is stochastically complete, or equivalently
	for every $\lambda>0$ the only non-negative bounded smooth solution
	$u$ of $\Delta u\geq\lambda u$ on $M$ is zero.
\end{theorem}
Note that the equivalence between stochastical completeness and the applicability of the
maximum principle is due to Grigor'Yan \cite{Gri99A}.

The following proposition will ensure that we can apply our maximum principle on $Z$:
\begin{proposition}\label{KhasTestMin}
	Let $Z\subset\X$ be a proper submanifold of $\X$ complete contractible space
	with sectional curvature pinched between $-b$ and $-a$, $b\geq a\geq 0$.
	Fix $p\in\X$ and denote $f$ the distance function to the point $p$.
	Then out of a compact set $K$, there is a constant $C>0$ such that the restriction of
	$f$ to $Z$ satisfies
	\begin{equation}
		\Delta(f|_Z)\leq C
	\end{equation}
\end{proposition}
\begin{proof}
	As~$\X$ is negatively pinched, there are~$b\geq a\geq 0$ such that~$-b\leq K_\X\leq -a$.
	By the Hessian comparison Theorem, the Hessian of $f$ is bounded between the hessian of
	the distance functions in the spaceforms of sectional curvature $-b$ and $-a$.
	Explicitely, this means
	\begin{equation}
		a\coth(af)g_\X\leq\nabla^2 f\leq b\coth(b f)g_\X
	\end{equation}
	If $a$ is zero, one has to replace $a\coth(af)$ by $\frac{1}{f}$, but it doesn't change
	anything to the proof.
	In particular, out of a compact set $K$ containing $p$,
	the hessian of $f$ is bounded by a constant $C$.
	Now because $Z$ is proper, the intersection $K\cap Z$ is a compact subset of $Z$.
	Also, because $Z$ is minimal, the laplacian of the restriction of $f$ satisfies
	\begin{equation}
		\Delta(f|_Z)=\Tr(\nabla^2 f|_{TZ})\leq C\dim Z\,.
	\end{equation}
	as desired.
\end{proof}
Applying theorem \ref{KhasMinskiiTest}, we deduce that the maximum principle is applicable
on a proper submanifold of such a space $\X$.
\begin{corollary}\label{WeakMaxPr}
	Let $\X$ be a complete contractible space whose sectional curvature is bounded
	between $-b$ and $0$ for some $b>0$.
	Let $Z\subset\X$ a proper minimal submanifold.
	Then for any $u\in C^2(Z)$ such that $\sup u<\infty$,
	for every $c<\sup u$,
	\begin{equation}
		\inf_{z\in Z: u(z)> c}\Delta u\leq 0
	\end{equation}
\end{corollary}
\begin{proof}
	Thanks to proposition \ref{KhasTestMin}, the function $f|_Z$ tends to $\infty$ at $\infty$,
	and out of a compact set $K$, it satisfies
	\begin{equation}
		\Delta f\leq C\leq C'f
	\end{equation}
	for some constant $C'>0$.
	In particular, we can apply theorem \ref{KhasMinskiiTest} to deduce that $Z$ is
	stochastically complete, which is equivalent to the statement of our corollary
	by Theorem~$3.1.$ of \cite{PRS05}.
\end{proof}
We now have everything needed to prove our main theorem.

\begin{proof}[Proof of theorem \ref{PlaPbAF}]
Introduce $u:Z\rightarrow\R$, such that $u(z)=d(z,Y)^2$.
Thanks to proposition \ref{BddDist}, $u$ is bounded.
But thanks to proposition \ref{SupHarmFct}, there is $C>0$ such that
\begin{equation}
\Delta u\geq C u
\end{equation}
	Now because of corollary \ref{WeakMaxPr},
	we can apply the maximum principle on $u$
	to get that $u=0$.
	Hence $Z\subset Y$, and by completeness $Z=Y$, as desired.
\end{proof}